\documentclass[12pt,a4paper]{article}

\usepackage[dvips]{graphics}
\usepackage[dvips]{epsfig}
\usepackage{latexsym}
\usepackage{amsfonts}
\usepackage{amssymb}
\usepackage{amsmath}
\usepackage[english]{babel}

\newtheorem{theorem}{Theorem}

\newtheorem{proposition}[theorem]{Proposition}
\newtheorem{definition}[theorem]{Definition}

\newtheorem{remark}[theorem]{Remark}

\newcommand{\dimo}{\noindent \emph{Proof. }}
\newcommand{\qed}{$\Box$\\}
\newcommand{\e}{\varepsilon}

\begin{document}

\title{\bf Cataloguing PL 4-manifolds \\ by gem-complexity}

\author{Maria Rita Casali\thanks{Work supported by the ``National Group for Algebraic and Geometric Structures, and their Applications" (GNSAGA - INDAM)
and by M.I.U.R. of Italy (project ``Strutture Geometriche, Combinatoria e loro Applicazioni'').}
\qquad Paola Cristofori\\
\small Dipartimento di Scienze Fisiche, Informatiche e  Matematiche\\[-0.8ex]
\small Universit\`a di Modena e Reggio Emilia\\[-0.8ex]
\small Modena, Italy\\
\small\tt casali@unimore.it \qquad paola.cristofori@unimore.it
}

\maketitle

\begin{abstract}
  We describe an algorithm to subdivide automatically a given set of PL $n$-manifolds (via {\it coloured triangulations} or, equivalently, via {\it crystallizations}) into classes whose elements are PL-homeomorphic.  The algorithm, implemented in the case $n=4$, succeeds to solve completely the PL-homeomorphism problem among the catalogue of all closed connected PL 4-manifolds up to gem-complexity 8 (i.e., which admit a coloured triangulation with at most 18 4-simplices).
\par
Possible interactions with the (not completely known) relationship among the different classifications in the TOP and DIFF=PL categories are also investigated. As a first consequence of the above PL classification, the non-existence of exotic PL 4-manifolds up to gem-complexity 8 is proved. Further applications of the tool are described, related to possible PL-recognition of different triangulations of the $K3$-surface.
\end{abstract}



\bigskip

\section{Introduction and main results}\label{introduction}

One of the most interesting features of piecewise-linear (PL) topology is the possibility of representing manifolds by combinatorial structures; the main developed theories concern the 3-dimensional case, where - thanks to recent advances in computing power - topologists succeeded in constructing exhaustive tables of ``small" 3-manifolds based on different representation methods (see \cite{[M_libro]} and its bibliography for successive results about closed orientable irreducible 3-manifolds  up to Matveev's complexity $11$, and \cite{[AM_2]}, \cite{[B]}  for analogous studies about closed non-orientable $\mathbb P^2$-irreducible 3-manifolds  up to Matveev's complexity $10$).

In dimension four, fewer combinatorial tools are available to represent PL-manifolds. On the other hand, classification results for topological (TOP) simply-connected $4$-mani\-folds are well-known, though more attention must be paid when considering equivalence of PL-structures.

{\it Crystallization theory} is a representation theory for PL-manifolds of arbitrary dimension by means of suitable edge-coloured graphs (called {\it crystallizations}), which are dual to coloured triangulations. Together with the Italian school that gave rise to the graph-theoretical tool  (see \cite{[P]}, \cite{[F]}, \cite{[G]}, \cite{[FG]}, \cite{survey} and their references), many authors around the world concurred to its development, with recent significant contributions, too: for example, \cite{[BaD]}, \cite{[BaS]}, \cite{[Sw]}.
The totally combinatorial nature of the representing objects and the generality with respect to dimension are among the strong points of crystallization theory: topological and PL properties are reflected into combinatorial ones, and the problem of distinguishing manifolds (both in TOP and in PL category) can be simplified by combinatorial invariants computed on the graphs.

In particular, in dimension four and five the best achievements have been obtained as regards the attempts of classifying PL-manifolds via a suitable graph-defined invariant, called {\it regular genus}\footnote{In the orientable case, it is the minimum genus of a surface where a graph representing the manifold regularly embeds: see \cite{[G]} for details.}:  they concern both the case of ``low" regular genus, and the case of ``restricted gap" between the regular genus of the manifold and the regular genus of its boundary, and the case of ``restricted gap" between the regular genus and the rank of the fundamental group of the manifold (see, for example, \cite{[CG]}, \cite{Casali} and \cite{Casali-Malagoli}).

More recently, the interest focused on other combinatorial invariants internal to crystallization theory, i.e. {\it GM-complexity} and {\it gem-complexity}, which are related in dimension three to Matveev's complexity, too (see \cite{[C$_4$]}, \cite{[CC$_1$]}, \cite{[CCM]}, \cite{[CC$_3$]}).
In particular, {\it gem-complexity} is the natural invariant used to create automatic catalogues of PL-manifolds via crystallizations: in fact, it is related to the minimum order of a crystallization of the manifold.
On the other hand, suitable moves on edge-coloured graphs are defined, which preserve the represented manifold up to PL-homeomorphisms; even if they are not able to solve algorithmically the recognition problem for general PL-manifolds, they are a powerful tool to face the problem itself, for a given set of PL-manifolds. In dimension three, this approach already allowed the development of a classification algorithm, which succeeded to completely recognize PL-homeomorphism classes of all $3$-manifolds up to gem-complexity $14$ (i.e. representable by coloured triangulations with at most $30$ tetrahedra): see \cite{[L]}, \cite{[CC$_1$]} and \cite{[CC$_2$]} for the orientable case and \cite{[C$_2$]}, \cite{[C$_4$]} and \cite{[BCrG$_1$]}  for the non-orientable one.

The present paper describes the $n$-dimensional extension of the above classifying algorithm, together with the results obtained by applying it to the crystallization catalogue representing all PL 4-manifolds up to gem-complexity $8$ (i.e. whose associated coloured triangulations have at most $18$ $4$-simplices).

The main classification results are collected into the following theorem, where $k(M^4)$ denotes the gem-complexity of $M^4$; the last statement makes use also of a partial analysis (whose completion is in progress) of the crystallization catalogue representing all PL $4$-manifolds with gem-complexity 9, which has been generated, too.

\begin{theorem}\label{Thm:1}
\ \
Let $M^4$ be a handle-free closed connected PL 4-manifold. Then:
\ \par
\begin{itemize}
\item{} $k(M^4)=0\ \Longleftrightarrow\ M^4$ is PL-homeomorphic to $\mathbb S^4$;
\item{} $k(M^4)=3\ \Longleftrightarrow\ M^4$ is PL-homeomorphic to $\mathbb{CP}^{2}$;
\item{} $k(M^4)=6\ \Longleftrightarrow\ M^4$ is PL-homeomorphic to either $\mathbb{S}^{2} \times \mathbb{S}^{2}$ or $\mathbb{CP}^{2} \# \mathbb{CP}^{2}$ or $\mathbb{CP}^{2} \# (- \mathbb{CP}^{2})$;
\item{} $k(M^4)=7\ \Longleftrightarrow\ M^4$ is PL-homeomorphic to $\mathbb{RP}^{4}$.
    \end{itemize}
Moreover:
\ \par
\begin{itemize}
\item{} no handle-free PL $4$-manifold $M^4$ exists with $k(M^4) \in \{1,2,4,5,8\};$
\item{} no exotic PL $4$-manifold exists, with $k(M^4) \le 8;$
\item{} any PL $4$-manifold $M^4$ with $k(M^4)=9$ is simply-connected (with second Betti number $\beta_2 \le 3$).
\end{itemize}
\end{theorem}

As far as the TOP category is concerned, the combinatorial properties of crystallizations, together with well-known results on TOP simply-connected 4-manifolds, yield the following interesting result related to the topological classification of simply-connected PL 4-manifolds with respect both to  {\it gem-complexity} and to  {\it regular genus}:

 \begin{theorem}\label{Thm:2}
 \label{Thm:classif_TOP general}  Let $M^4$ be a simply-connected PL 4-manifold $M^4$.
 If either gem-complexity $k(M^4)\leq 65$ or regular genus $\mathcal G(M^4)\leq 43$, then $M^4$ is TOP-homeomorphic to
 $$(\#_r\mathbb {CP}^2)\#(\#_{r^\prime}(-\mathbb {CP}^2))\quad or\quad\#_s(\mathbb S^2\times \mathbb S^2),$$
 \noindent where $r+r^\prime = \beta_2(M^4),\ s = \frac 1 2\beta_2(M^4)$ and $\beta_2(M^4)$  is the second Betti number of $M^4$.
 \end{theorem}

Theorem~\ref{Thm:2} summarizes Proposition~\ref{Thm:classif_TOP via gem-complexity} (for gem-complexity) and Proposition~\ref{Thm:classif_TOP via regular-genus} (for regular genus) of Section~\ref{sec:TOP-classification}.

\bigskip

As it is well-known, up to now there is no classification of smooth structures on any given smoothable topological 4-manifold; on the other hand, finding non-diffeomorphic smooth structures on the same closed simply-connected topological manifold has long been an interesting problem.

We hope that further advances in the generation and classification of crystallization catalogues for PL $4$-manifolds, according to gem-complexity, could  produce examples of non-equivalent PL-structures on the same topological 4-manifold. For example, if at least one among the infinitely many PL 4-manifolds TOP-homeomorphic but not PL-homeomorphic to $\mathbb{CP}^{2}\#_2(-\mathbb{CP}^{2})$  (which are proved to exist in  \cite{[A-DP$_2010$]}) admits a so called {\it simple crystallization}  (according to \cite{[BaS]}), then it will appear in the catalogue of order $20$ crystallizations.

Moreover, we point out that the program performing automatic recognition of PL-homeo\-morphic 4-manifolds may be a useful tool to approach open problems related to different triangulations of the same TOP 4-manifold, which are conjectured to represent the same PL 4-manifold, too. The first candidates are the two known 16-vertices and 17-vertices triangulations of the $K3$-surface: see \cite{[CK]}  and \cite{[SK]}, together with the attempts to settle the conjecture described in \cite{[BuS]}, \cite{[BL]}  and \cite{[BaS]}.

\section{Basic notions of crystallization theory}

As already pointed out, \textit{crystallization theory} allows to represent combinatorially PL manifolds of arbitrary dimension, without restrictions concerning orientability, connectedness or boundary properties, by means of suitable \textit{edge-coloured graphs} or - equivalently - by means of \textit{coloured triangulations}. A detailed account of the theory may be found in  \cite{[FGG]},  \cite{[L]} and \cite{survey}, together with their references.

\medskip

In the present paper, when not otherwise stated, we will restrict our attention to the case of closed, connected PL $n$-manifolds.

\begin{definition}
An \textit{(n+1)-coloured graph }
is a pair $(\Gamma,\gamma),$ where $\Gamma=(V(\Gamma), E(\Gamma))$
is a regular multigraph\footnote{According to \cite{[W]}, this means that all vertices of $V(\Gamma)$ have the same degree, that loops are forbidden, while multiple edges are allowed.} of degree $n+1$ and
$\gamma : E(\Gamma) \to \Delta_n=\{0,1,\dots,n\}$ is injective on
adjacent edges.
\end{definition}

The elements of the set $\Delta_n=\{0,1,\dots,n\}$ are called
\textit{colours}; moreover, for each $i\in\Delta_n$, we denote by
$\Gamma_{\hat i}$ the $n$-coloured graph obtained from
$(\Gamma,\gamma)$ by deleting all edges coloured by
$i$.

\begin{definition} An $(n+1)$-coloured graph  $(\Gamma,\gamma)$ is said to be {\it contracted} if the subgraph $\Gamma_{\hat i}$ is connected, for each $i \in \Delta_n.$
\end{definition}

Each $(n+1)$-coloured graph uniquely determines an $n$-dimensional CW-complex $K(\Gamma)$, which is said to be {\it associated to $\Gamma$}:
\begin{itemize}
\item{} for every vertex $v\in V(\Gamma)$, take an $n$-ball $\sigma(v)$  abstractly isomorphic to an $n$-simplex, and label injectively its $n+1$ vertices by the colours of $\Delta_n$;
\item{} for every $i$-coloured edge between $v,w\in V(\Gamma)$, identify the ($n-1$)-faces of $\sigma(v)$ and $\sigma(w)$ opposite to $i$-labelled vertices, so that equally labelled vertices
coincide.
\end{itemize}

\bigskip

It is easy to check that the properties of $(\Gamma,\gamma)$ imply $K(\Gamma)$ to be actually a {\it pseudo-complex}: in particular, its balls may intersect in more than one face, but self-identifications of faces are not allowed.\footnote{Even if, in general, $K(\Gamma)$ fails to be a simplicial complex, we will always call {\it $h$-simplices} its $h$-balls, for every $h\le n$.}

\begin{definition} An $(n+1)$-coloured graph $(\Gamma,\gamma)$ is said to {\it represent} a PL $n$-manifold $M^n$ (briefly, it is a {\it gem} of $M^n$) if $M^n$ is PL-homeomorphic to $|K(\Gamma)|$.  If, in addition, $(\Gamma,\gamma)$ is contracted, then it is called a {\it crystallization} of $M^n$. In both cases, the pseudo-complex $K(\Gamma),$ equipped with the vertex-labelling inherited from $\gamma$, is said to be
a {\it coloured triangulation} of $M^n$.\end{definition}

\begin{remark} \label{rem:baricentrica} It is easy to prove that each PL $n$-manifold admits an $(n+1)$-coloured graph representing it: just take the barycentric subdivision
$H^{\prime}$ of any (simplicial) triangulation $H$ of $M^n$, label any vertex of $H^{\prime}$ with the dimension of the open simplex containing it in $H$, and then consider the 1-skeleton of the dual cellular complex of $H^{\prime}$, with each edge coloured by $i$ iff it is dual to an $(n-1)$-face whose vertices are labelled by $\Delta_n-\{i\}$.
\end{remark}

The following proposition collects some well-known facts in crystallization theory:

\begin{proposition} \label{Thm:raccolta_preliminari}
Let $(\Gamma,\gamma)$ be an $(n+1)$-coloured graph.
\begin{itemize}
\item[(a)] $|K(\Gamma)|$ is orientable iff  $\Gamma$ is bipartite.
\item[(b)] $\forall \mathcal B \subset \Delta_n$, with $ \# \mathcal B=h,$
there is a bijection between ($n-h$)-simplices of $K(\Gamma)$
whose vertices are labelled by $\Delta_n - \{\mathcal B\}$ and
connected components of the
$h$-coloured graph $\Gamma_{\mathcal
B}=(V(\Gamma),$ $\gamma^{-1}(\mathcal B))$  (which are called \emph{$h$-residues involving colours $\mathcal B$}, or \emph{$\mathcal B$-residues} of $\Gamma$, and whose number will be denoted by $g_{\mathcal B}$).
\par In particular: \ $c$-labelled vertices of
$K(\Gamma)$ are in bijection with connected components of $\Gamma_{\hat c}= \Gamma_{\Delta_n -\{c\}}.$
\item[(c)] If $\#V(\Gamma) = 2p$ and $\sum_{\#\mathcal B=h} g_{\mathcal B}$ denotes the total number of $h$-residues of $\Gamma$,
    then $$\chi (|K(\Gamma)|) =
        (-1)^{n-1} \cdot p \cdot (n-1) + \sum_{h=2}^n (-1)^h  \cdot \sum_{\#\mathcal B=h}  g_{\mathcal B}.$$
\item[(d)] $|K(\Gamma)|$ is an $n$-manifold if and only
if, for every $c \in \Delta_n$, each connected component of $
\Gamma_{\hat c}$ represents $\mathbb S^{n-1}$.
\item[(e)] If $(\Gamma,\gamma)$ is a crystallization of an $n$-manifold $M^n$, then $$rk (\pi_1(M^n)) \le \min \{g_{\mathcal B} -1  \ | \ \# \mathcal B = n-1\}.$$
\end{itemize}
\end{proposition}

\par \noindent Moreover:
\begin{proposition} [{\bf Pezzana Theorem}]
Each PL $n$-manifold admits a crystallization.
\end{proposition}

It is not difficult to understand that, generally, many crystallizations of the same PL $n$-manifold exist; hence, it is a basic problem how to recognize crystallizations (or, more generally, gems) of the same PL $n$-manifold.

The easiest case is that of two \textit{colour-isomorphic} gems,
i.e. if there exists an isomorphism between the graphs, which
preserves colours up to a permutation of $\Delta_n$. It is quite
trivial to check that two colour-isomorphic gems produce the same
polyhedron.

The following result assures that colour-isomorphic graphs can be
effectively detected by means of a suitably defined numerical
\textit{code}, which can be directly computed on each of them (see
\cite{[CG$_2$]}).

\begin{proposition}
Two $(n+1)$-coloured graphs are colour-isomorphic iff their codes coincide.
\end{proposition}

\begin{figure}[th]
\begin{center}\scalebox{0.45}{\includegraphics{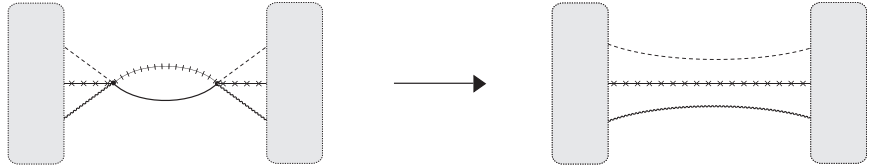}}\end{center}
\caption{\label{fig:dipole_move} dipole move}
\end{figure}

The problem of recognizing non-colour-isomorphic gems
representing the same manifold is also solved, but not algorithmically:
a finite set of moves - the so called \textit{dipole moves} - is
proved to exist, with the property that two gems represent the
same manifold iff they can be related by a finite sequence of such
moves. A dipole move consists in the insertion or elimination of particular configurations involving $h$ parallel edges, called {\it $h$-dipoles} ($1 \le h \le n$): see \cite{[FG]} for details, or Figure~\ref{fig:dipole_move}
for an example in dimension $n=4$, with $h=2$.

\medskip

\bigskip

In this paper, however, we will also make use of another set of moves, which appears to be more suitable for algorithmic procedures (see the notion of {\it blob} and {\it flip} in Section~\ref{sec:paragrafo algoritmo classificazione}). Even if they still do not solve algorithmically the problem for general PL $n$-manifolds (nor for general PL  $4$-manifolds), nevertheless we will prove that a fixed sequence of  blobs and flips is sufficient to classify - via PL-homeomorphism - all PL 4-manifolds admitting a coloured triangulation with at most 18 $4$-simplices (see Section~\ref{sec:PL classification}).

\bigskip

In order to define the class of gems involved in our catalogues, further preliminary notions are required.

\begin{definition} A pair $(e,f)$ of equally coloured edges in an $(n+1)$-coloured graph $(\Gamma,\gamma)$ is said
to form a \textit{$\rho_s$-pair} iff $e$ and $f$ both belong to exactly $s$ common bicoloured cycles of $\Gamma$.
\end{definition}

\begin{figure}[th]
\begin{center}\scalebox{0.45}{\includegraphics{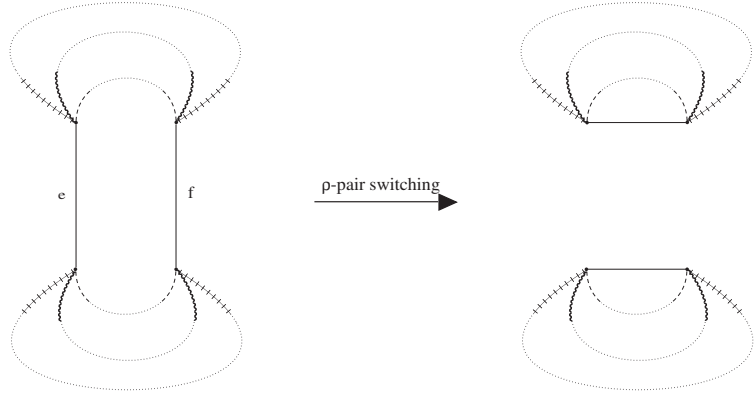}}\end{center}
\caption{\label{fig:rho-switching} $\rho$-pair switching}
\end{figure}

Figure~\ref{fig:rho-switching} shows the combinatorial move called {\it $\rho$-pair switching}.

The effect of $\rho$-pair switching on crystallizations is explained by the following result, where $\mathbb H$ denotes an {\it $n$-dimensional handle}, i.e. either the orientable or non-orientable $\mathbb S^{n-1}$-bundle over $\mathbb S^1$ (respectively denoted by $\mathbb{S}^{1} \times \mathbb{S}^{n-1} $ and  $\mathbb{S}^{1} \widetilde \times \mathbb{S}^{n-1}$), according to the orientability of $M^n$:

 \begin{proposition} \label{Thm:rho-pairs} {\rm (\cite{[BG_2011]})} Let $(\Gamma,\gamma)$ be a crystallization of a PL $n$-manifold $M^n$, $n\geqslant 3$ and let $(\Gamma^{\prime}, \gamma^{\prime})$ be obtained by switching a $\rho_s$-pair in $\Gamma$. Then:
\begin{itemize}
\item[(a)] if $s=n-1,$ \ $(\Gamma^{\prime}, \gamma^{\prime})$ is a gem of $M^n$, too;
\item[(b)] if $s=n,$ \  $(\Gamma^{\prime}, \gamma^{\prime})$ is a gem of an $n$-manifold $N^{n}$ such that $M^n \cong_{PL} N^{n} \# \mathbb H.$
\end{itemize}
\end{proposition}

\begin{definition} An $(n+1)$-coloured graph is said to be \textit{rigid} (resp. \textit{rigid dipole-free}) if it has no  $\rho_s$-pairs, with $s \in \{n-1,n\}$ (resp. if it is rigid and has no dipoles).
\end{definition}

Catalogues of PL $n$-manifolds are obviously constructed with respect to increasing ``complexity" of the representing combinatorial objects. Within crystallization theory, the following quite natural invariant is considered:

\begin{definition} Given a PL $n$-manifold $M^n$, its \textit{gem-complexity} is the non-negative integer $k(M^n)= p - 1$, where $2p$ is the minimum order (i.e. the minimum number of vertices) of a crystallization of $M^n$.
\end{definition}

In order to generate an exhaustive catalogue of $(n+1)$-coloured graphs representing all closed connected PL $n$-manifolds up to a fixed gem-complexity, the restriction to the class of rigid dipole-free crystallizations yields no loss of generality, as Proposition~\ref{Thm:rigid&dipole-free} below proves.

\medskip

In the following,  $\#_h M$ denotes the connected sum of $h$ copies of a given $n$-manifold $M$.

\begin{definition} A PL $n$-manifold $M^n$ is said to be {\it handle-free} if it admits no handles as connected summands.
\end{definition}

\begin{proposition} \label{Thm:rigid&dipole-free}
Let $M^n$ be a PL $n$-manifold ($n\ge 3$).
Then:
\begin{itemize}
\item[(a)] If $M^n$ is handle-free, then a rigid dipole-free order $2p$ crystallization $(\Gamma, \gamma)$ of $M^n$ exists, so that $k(M^n)=p-1$.
\item[(b)] Otherwise, a  rigid dipole-free order $2p$ crystallization of a PL $n$-manifolds $N^n$ exists, so that $M^n \cong_{PL} N^n \#_h \mathbb H$ ($h \geq 0$) and $k(M^n)= p - 1 + n \cdot h.$
  \end{itemize}
\end{proposition}

\dimo
Statement (a) directly follows from \cite[Theorem~5.3]{[BG_2011]}.

Let now $M^n \cong_{PL} N^n \#_h \mathbb H$ be any connected sum decomposition of $M^n$, with $h \geq 0$.\footnote{Note that no assumption is made, both on $h$ and on $N^n$: the decomposition turns out to be ``trivial" in case $h=0$ and $N^n = M^n$. Note also that in dimension $n=3$, where the uniqueness of the decomposition under connected sum holds, the statement is already known: see  \cite[Proposition~8(b)]{[C$_2$]}.}
Since a standard order $2(n+1)$ crystallization of $\mathbb{S}^{1}\times\mathbb{S}^{n-1}$ (resp. $\mathbb{S}^{1} \widetilde \times \mathbb{S}^{n-1}$) is well-known (see \cite{[GV]}), and since the so called {\it graph connected sum} (see \cite{[FGG]}) yields an order $2(p_1+p_2-1)$  gem of $N_1\#N_2$ from any order $2p_1$ (resp. $2p_2$) gem of $N_1$ (resp. $N_2$),  inequality $k(M^n) \leq p-1 + n \cdot h$ easily follows, $2p$ being the order of any crystallization of $N^n.$

Hence, if we set
$$ \begin{aligned} \mathcal S = \Big\{ p - 1 + n \cdot h \ \ | & \ \ M^n \cong_{PL} N^n \#_h \mathbb H, \ \ 2p= \# V(\Gamma), \ \\ \ & \ \ \Gamma \ \text{rigid dipole-free crystallization of} \ N^n \Big\}, \end{aligned}$$
it is proved that
$$k(M^n) \, \leq \, \min \mathcal S,$$
where the minimum is taken over all decompositions $M^n \cong_{PL} N^n \#_h \mathbb H$ (with $h\geq0$) and over all rigid dipole-free crystallizations of $N^n.$

In order to prove the reversed inequality (and hence statement (b)), let us consider an arbitrary crystallization $(\bar \Gamma, \bar \gamma)$ of $M^n$, with order $2 \bar p$.  If $(\bar \Gamma, \bar \gamma)$ is a rigid dipole-free crystallization, then  $\bar p -1 \in \mathcal S$ (with respect to the trivial decomposition of $M^n$), and so  $\bar p -1 \geq  \min \mathcal S$ trivially holds.
If $(\bar \Gamma, \bar \gamma)$ contains $\rho_{n-1}$-pairs and/or dipoles, a order $2 {\bar p}^{\prime}$ rigid dipole-free crystallization $(\bar \Gamma^\prime, \bar \gamma^\prime)$ of $M^n$ is easily obtained via a suitable number of $\rho_{n-1}$-switching, each one followed by a 1-dipole elimination, and/or dipole eliminations; hence, $ {\bar p}^{\prime} -1 \in \mathcal S$ (with respect to the trivial decomposition of $M^n$), and so $\bar p -1 > {\bar p}^{\prime} -1 \geq \min \mathcal S$ trivially holds.

Finally, let us assume $(\bar \Gamma, \bar \gamma)$ to admit no $\rho_{n-1}$-pairs, no dipoles and $h\geq 1$ $\rho_n$-pairs.
After each $\rho_n$-pair switching, $n$ 1-dipoles appear, one for each colour not involved in the $\rho$-pair. Let $(\bar \Gamma^{\prime\prime}, \bar \gamma^{\prime\prime})$ be the order $2 {\bar p}^{\prime\prime}$ graph obtained by all $\rho_n$-switchings and elimination of the resulting $n \cdot h$ 1-dipoles;
according to Proposition~\ref{Thm:rho-pairs}(b), $(\bar \Gamma^{\prime\prime}, \bar \gamma^{\prime\prime})$ is a rigid dipole-free crystallization of a PL $n$-manifold ${\bar N}^n$ such that $M^n \cong_{PL} {\bar N}^n \#_h \mathbb H.$ Hence, $ {\bar p}^{\prime\prime} -1 + n \cdot h \in \mathcal S.$
Relation $\bar p -1 = {\bar p}^{\prime\prime} -1 + n \cdot h \geq \min \mathcal S$ easily follows.\hskip 340pt\qed

Elementary notions of crystallization theory - and in particular the graph-connected sum quoted in the above proof - allow to easily extend to any dimension both the {\it sphere-recognition property} and the {\it finiteness property} and the {\it sub-additivity (with respect to connected sum)} of gem-complexity, already stated in \cite[Proposition~10]{[C$_2$]} for the $3$-dimensional case.

\medskip

Although the gem-complexity cannot be additive on the whole set of PL 4-manifolds, as a consequence of Wall theorem (\cite{Wall}), nevertheless additivity is possible on restricted classes of manifolds: for example those with small values of gem-complexity ($\leq 8$ as shown by the results of the present paper) or those admitting {\it simple} crystallizations (as proved in \cite{[CCG]}).

The invariant regular genus shares with the gem-complexity the properties of sphere recognition and subadditivity in any dimension. In dimension three it is additive but not finite-to-one (since it coincides with the Heegaard genus), while in the closed 4-dimensional case both finiteness and additivity are open problems. However, additivity holds within the class of 4-manifolds admitting {\it simple} crystallizations (\cite{[CCG]}).

\begin{remark} Note that, as pointed out in \cite[Remark~1]{[FG$_2$]}, additivity of the regular genus in dimension four would imply the Smooth Poincar\'e Conjecture.
\end{remark}

\section{4-dimensional generation algorithm} \label{sec:paragrafo algoritmo generazione}

By Proposition~\ref{Thm:raccolta_preliminari}(d), the generation of catalogues of crystallizations of $n$-manifolds with a fixed number of vertices $2p$ is essentially inductive on dimension and requires the prior generation and recognition of all gems (with $2p$ vertices) representing the $(n-1)$-sphere.

It is therefore easy to understand why the first results have been obtained in dimension three, since 2-sphere recognition can be performed easily by computing the Euler characteristic, and null Euler characteristic characterizes closed 3-manifolds.

However, the generating algorithm, even in low dimension, becomes quickly very intensive as the number of vertices increases and requires large computing resources.
A way to face this problem is to find combinatorial configurations in the graphs, which can be eliminated without changing the manifold. Examples of such configurations are dipoles and $\rho$-pairs.

Proposition~\ref{Thm:rigid&dipole-free} assures that restricting the catalogues to rigid dipole-free crystallizations does not affect their completeness.

Catalogues of 3-manifold rigid\footnote{It is easy to see that in dimension three, rigidity and contractedness imply absence of dipoles.} crystallizations up to 32 vertices have already been generated and completely classified (\cite{[C$_2$]}, \cite{[CC$_1$]}, \cite{[CC$_2$]}, \cite{[BCrG$_1$]}). Complete classification was also obtained for genus two 3-manifolds admitting a crystallization with at most 42 vertices (\cite{[BCrG$_2$]}).

\smallskip
Let now fix our attention to the 4-dimensional case.
For each $p \ge 1$, we will denote by $\mathcal C^{(2p)}$ (resp. $\tilde {\mathcal C}^{(2p)}$) the catalogue of all not colour-isomorphic rigid dipole-free bipartite (resp. non-bipartite) crystallizations of 4-manifolds with $2p$ vertices.

Note that if $\Gamma\in\mathcal C^{(2p)}\cup\tilde {\mathcal C}^{(2p)}$, then $\Gamma_{\hat 4}$ is a (not necessarily contracted, nor rigid) 4-coloured graph representing $\mathbb S^3$ and lacking in $\rho_3$-pairs. Let $S^{(2p)}$ denotes the set of such 4-coloured graphs.

$S^{(2p)}$ will be the starting set of the procedure generating $\mathcal C^{(2p)}$ and $\tilde {\mathcal C}^{(2p)}$, which consists essentially in adding 4-coloured edges to all elements of $S^{(2p)}$, so as to obtain crystallizations of 4-manifolds.

The set $S^{(2p)}$ is constructed by a suitable adaptation of the 3-dimensional generation algorithm; recognition of the 3-sphere is performed by cancelling dipoles and switching $\rho$-pairs in order to obtain a rigid crystallization and by comparing the resulting graph with the list of rigid crystallizations of $\mathbb S^3$, which appear in the 3-dimensional catalogue.

As a matter of fact, the recognition is very easy for $p<12$, since the only rigid crystallization of $\mathbb S^3$ up to this order is the standard one with two vertices.

So, the generating algorithm in dimension four runs as follows:

\begin{itemize}\item[1)]
Construct the set $S^{(2p)}=\{ \Sigma^{(2p)}_1,
\Sigma^{(2p)}_2, \dots, \Sigma^{(2p)}_{n_p}\}$.
\item[2)] For each $i=1,2,\dots,n_p:$
\begin{itemize}
\item[-] add to $\Sigma^{(2p)}_i$ 4-coloured edges in all possible ways so to produce 4-coloured graphs;
\item [-] for each produced graph $\Gamma$, check absence of $\rho$-pairs and 2-dipoles;
\item [-] for each $c\in\Delta_3$, check that $\Gamma_{\hat c}$ represents $\mathbb S^3$.
\end{itemize}
\item [3)] Compute and compare the codes in order to exclude colour-isomorphic duplicates.
\end{itemize}

However, the above algorithm is practically useless due to the great computational time it requires; therefore it needs some modifications in order to be effective. More precisely, a branch and bound technique is used to prune the tree of possible attachments of edges on each element of $S^{(2p)}$.

Let $\bar\Gamma$ be a 5-coloured graph obtained from an element of $\mathcal S^{(2p)}$ by addition of $r<p$ 4-coloured edges, then $\bar\Gamma$ will be kept for further additions if and only if:

\begin{itemize}
\item [(i)] it contains no three edges with the same endpoints (otherwise there will be $\rho$-pairs in the final regular graph);
\item [(ii)] for each $i\in\Delta_3$, $\bar\Gamma_{\hat i}$ represents a $3$-sphere with holes.
\end{itemize}

Unfortunately condition (ii) is very heavy to check, since it implies recognition of 3-spheres with holes.  Instead, we use a weaker condition, which is equivalent to requiring $\bar\Gamma_{\hat i}$ to be a $3$-manifold (with boundary), i.e.

\begin{itemize}
\item [(ii$^\prime$)] for each pair of colours $i,j\in\Delta_3$, each 3-residue of $\bar\Gamma$ not involving colours $i, j$ must represent a disjoint union of 2-spheres, possibly with holes.
\end{itemize}

Condition (ii$^\prime$) can be checked by direct computation on $\bar\Gamma$ in the following way.

Note that $\bar\Gamma$ is a so-called {\it 5-coloured graph with boundary}, i.e. $\bar\Gamma_{\hat 4}$ is a 4-coloured graph. Then $K(\bar\Gamma)$, which is obtained exactly in the same way as in the closed case, is a pseudocomplex with non-empty boundary. The 3-simplices which triangulate $\partial K(\bar\Gamma)$ correspond bijectively to the vertices of $\bar\Gamma$ missing the 4-coloured edge ({\it boundary vertices}).

A 4-coloured graph (without boundary) $\partial\bar\Gamma$ such that $|K(\partial\bar\Gamma)|\cong |\partial K(\bar\Gamma)|$ can be constructed as follows:
\begin{itemize}
\item $V(\partial\bar\Gamma)$ is the set of boundary vertices of $\bar\Gamma$;
\item for each $c\in\Delta_3$, two vertices of $\partial\bar\Gamma$ are $c$-adjacent iff they are connected by a $\{c,4\}$-coloured path in $\bar\Gamma$.
\end{itemize}

Then, a suitable extension of the Euler characteristic computation of Proposition~\ref{Thm:raccolta_preliminari}(c) implies that condition (ii') is equivalent to requiring the following equality to hold:

\begin{equation}\label{eq:planarity4}
\sum_{k,t\in\Delta_4\setminus\{i,j\}} g_{kt}-\frac{\bar p}2 = 2g_{\hat i \hat j} - \bar g_{\hat i \hat j}\ ,
\end{equation}

\noindent where $ g_{kt}$ is the number of $\{k,t\}$-coloured cycles of $\bar\Gamma$, $\bar p$ is the number of boundary vertices of $\bar\Gamma$ and $g_{\hat i \hat j}$ (resp. $\bar g_{\hat i \hat j}$) is the number of $(\Delta_4 -\{i,j\}$)-residues of $\Gamma$ (resp. $\partial \Gamma$).

The above described restrictions succeed in reducing considerably both the computation time and the size of the resulting catalogues. Moreover, a parallelization strategy, which has been adopted in the implementation, has allowed to reduce further the computation time: see \cite{[MRC]} for details.

\smallskip
As a consequence we could produce catalogues $\mathcal C^{(2p)}$ and $\tilde {\mathcal C}^{(2p)}$ for each $p\leq 10$ (see Table~1 below).

\bigskip

\centerline{\begin{tabular}{|c|| c | c | c | c | c | c | c | c | c |  c |}
 \hline
  {\scriptsize \bf 2p } & \scriptsize{2} & \scriptsize{4} & \scriptsize{6} & \scriptsize{8} & \scriptsize{10} & \scriptsize{12} & \scriptsize{14} & \scriptsize{16} & \scriptsize{18}  & \scriptsize{20}
  \\   \hline  \hline  \ & \ & \ & \ & \ & \ & \ & \ & \ & \  & \
  \\
 {\scriptsize \bf  $\# S^{(2p)}$} \hfill & \scriptsize{1} & \scriptsize{0} & \scriptsize{2} & \scriptsize{9} & \scriptsize{39} & \scriptsize{400} & \scriptsize{5.255} & \scriptsize{95.870} & \scriptsize{1.994.962} & \scriptsize{45.654.630}
  \\  \hline  \hfill \ & \ & \ & \ & \ & \ & \ & \ & \ & \  & \
  \\
 \ & \ & \ & \ & \ & \ & \ & \ & \ & \  & \
  \\
{ {\scriptsize \bf  $\#\ \mathcal C^{(2p)}$}}
 \hfill & {\scriptsize{1}} & {\scriptsize{0}} & {\scriptsize{0}} & {\scriptsize{1}} & {\scriptsize{0}} & {\scriptsize{0}} & {\scriptsize{1.109}} & {\scriptsize{4.511}} & {\scriptsize{44.803}}  & { \scriptsize{47.623.129}}
\\  \hline  \hfill \ & \ & \ & \ & \ & \ & \ & \ & \ & \  & \
  \\
 \ & \ & \ & \ & \ & \ & \ & \ & \ & \  & \
  \\
{ {\scriptsize \bf  $\#\ \mathcal{\tilde C}^{(2p)}$}}
 \hfill & {\scriptsize{0}} & {\scriptsize{0}} & {\scriptsize{0}} & {\scriptsize{0}} & {\scriptsize{0}} & {\scriptsize{0}} & {\scriptsize{0}} & {\scriptsize{1}} &  { \scriptsize{0}}  & { \scriptsize{0}}
 \\ \hline
  \end{tabular}}

\bigskip

\centerline{Table~1}

\begin{remark} We point out that the unique rigid dipole-free crystallization of $\mathcal C^{(2)}$ (resp. of $\mathcal C^{(8)}$) is the standard crystallization of $\mathbb S^4$ (resp. $\mathbb {CP}^2$: see \cite{[G$_1$]}), while the unique non-bipartite rigid dipole-free crystallization appearing up to $20$ vertices is the standard one of $\mathbb {RP}^4$ with $16$ vertices (\cite{[G$_2$]}).
\end{remark}

\begin{remark} A further restriction on the catalogues could be imposed: if $\Gamma\in\mathcal C^{(2p)}\cup\tilde {\mathcal C}^{(2p)}$ is a {\it graph connected sum} (see \cite{[FGG]}, or the proof of Proposition~\ref{Thm:rigid&dipole-free}) of two graphs $\Gamma_1$ and $\Gamma_2$, then $|K(\Gamma)| \cong_{PL} |K(\Gamma_1)|\# |K(\Gamma_2)|$ and we call $\Gamma$ {\it splittable}. The problem of recognizing  the manifold $|K(\Gamma)|$ is thus traced back to the (easier) problem of recognizing $|K(\Gamma_1)|$ and $|K(\Gamma_2)|$. However, the low number of splittable crystallizations (about 0.6\% of the elements for the catalogues with up to 18 vertices) forces to directly recognize the crystallizations by moves keeping the PL-homeomorphism of $|K(\Gamma)|$ (see Section~\ref{sec:paragrafo algoritmo classificazione}).
\end{remark}

\section{TOP classification via combinatorial invariants}  \label{sec:TOP-classification}

Before describing the classification performed in our crystallization catalogues via suitable sequences of moves which realize PL-homeomorphisms of the represented PL 4-manifolds, we devote the present section to the much weaker problem of classifying the involved PL 4-manifolds within the TOP category.

\bigskip

The starting point is the direct computation of the Betti numbers and of the rank of the fundamental group for each PL 4-manifold represented by a crystallization of our catalogues.

\begin{proposition}
\label{Thm:calcolo_Eulero}
All orientable PL 4-manifolds represented by elements of $\mathcal C^{(2p)}$,  $1 \le p \le 10,$ are simply-connected.
The unique elements of $\mathcal C^{(2)}$ and $\mathcal C^{(8)}$ are known crystallizations representing $\mathbb S^4$ and $\mathbb {CP}^2$ respectively.
\ \par \noindent
Moreover:
\begin{itemize}
\item[(a)] Among the $1.109$ crystallizations of $\mathcal C^{(14)}:$
\begin{itemize}
\item{} exactly one represents a (simply-connected) PL 4-manifold $M^4$ with $\beta_2=1;$
\item{} all the remaining ones represent (simply-connected) PL 4-manifolds $M^4$ with $\beta_2=2.$
\end{itemize}
\item[(b)] All $4.511$ crystallizations of $\mathcal C^{(16)}$ represent (simply-connected) PL 4-manifolds $M^4$ with $\beta_2=2.$
\item[(c)] Among the $44.803$ crystallizations of $\mathcal C^{(18)}:$
\begin{itemize}
\item{} ten represent (simply-connected) PL 4-manifolds $M^4$ with $\beta_2=1;$
\item{} all the remaining ones represent (simply-connected) PL 4-manifolds $M^4$ with $\beta_2=2.$
\end{itemize}
\item[(d)] Among the $47.623.129$ crystallizations of $\mathcal C^{(20)}:$
\begin{itemize}
\item{} exactly one represents a (simply-connected) PL 4-manifold $M^4$ with $\beta_2=0;$
\item{} $370$ represent (simply-connected) PL 4-manifolds $M^4$ with $\beta_2=1;$
\item{} $501.900$ represent (simply-connected) PL 4-manifolds $M^4$ with $\beta_2=2;$
\item{} all the remaining ones represent (simply-connected) PL 4-manifolds $M^4$ with $\beta_2=3.$
\end{itemize}
\end{itemize}
\end{proposition}

\dimo
In virtue of Proposition~\ref{Thm:raccolta_preliminari}(e), an estimation of  $rk(\pi_1(M^4))$, for each PL 4-manifold $M^4$ represented by our catalogues, may be obtained by computing the number $g_{rst}$ of  3-residues involving colours $\{r,s,t\}$, with $0 \le r < s < t \le 4$,  for each element  of $\mathcal C^{(2p)}$,  $1 \le p \le 10.$
The calculation has been done by means of a suitable procedure of the program DUKE III \footnote{{\it ``DUKE III: A program to handle edge-coloured graphs representing PL n-dimensional manifolds"} is available on the Web: http://cdm.unimo.it/home/matematica/casali.mariarita/DUKEIII.htm}
and the program output ensures that each crystallization of $\mathcal C^{(2p)}$,  $1 \le p \le 10,$ has $g_{rst}=1$ for at least a choice of distinct $r,s,t\in  \Delta_4.$
Hence, the simply-connectedness of all involved orientable PL 4-manifolds is proved.

\smallskip

In order to calculate the Betti numbers of the same PL 4-manifolds, it is necessary to apply Proposition~\ref{Thm:raccolta_preliminari}(c), yielding the Euler characteristic of $M^4$ by a direct computation on each order $2p$ crystallization of $M^4$:\footnote{Recall that contractedness implies $\Gamma$ to have exactly one $4$-residue involving $\Delta_4 - \{i\}$, for each $i \in \Delta_4$.}
\begin{equation}\label{eq:Euler}
\chi (M^4) =  5 - \sum_{i<j<k} g_{ijk} + \sum_{i<j} g_{ij} -  3p.
\end{equation}
Now, the simply-connectedness implies $ \chi (M^4) = 2 + \beta_2(M^4)$; hence, $\beta_2(M^4)$ follows by a direct computation of both the number of 3-residues and 2-residues, for each element  of $\mathcal C^{(2p)}$,  $1 \le p \le 10.$
The statement is proved by making use of suitable procedures of the program DUKE III.
\hskip 447pt\qed

The fact that all orientable PL 4-manifolds represented by crystallizations of our catalogues are simply-connected has important consequences as regards their classification in the TOP category.

In fact, the following result proves that, up to a significantly high gem-complexity, the classification of PL 4-manifolds up to TOP-homeomorphism is quite easy, at least in the simply-connected case (which - as a matter of fact - turns out to be the most frequent case):\footnote{The statement of Proposition~\ref{Thm:classif_TOP via gem-complexity} was announced in \cite{oberwolfach} and in \cite{[CC_EN]}.}

\begin{proposition}
\label{Thm:classif_TOP via gem-complexity} Any simply-connected PL 4-manifold $M^4$, with $k(M^4)\leq 65$, is TOP-homeo\-morphic to
$$(\#_r\mathbb {CP}^2)\#(\#_{r^\prime}(-\mathbb {CP}^2))\quad or\quad\#_s(\mathbb S^2\times \mathbb S^2),$$
\noindent where $r+r^\prime = \beta_2(M^4),\ s = \frac 1 2\beta_2(M^4)$ and $\beta_2(M^4) \le \frac{k(M^4)}3$ is the second Betti number of $M^4$.
\end{proposition}

\dimo
Let $\Gamma$ be an order $2p$ crystallization of $M^4$. As already recalled, formula (\ref{eq:Euler}) yields the direct computation of the Euler characteristic of $M^4.$

On the other hand, the planarity of each 3-residue of $\Gamma$ yields (via Proposition~\ref{Thm:raccolta_preliminari}(c), too)
$2 g_{ijk}  =  g_{ij} + g_{ik} + g_{jk} - p$ for each triple $(i,j,k)\in \Delta_4$, from which the following relation is obtained:
$$ 2 \sum_{i<j<k} g_{ijk} = 3 \sum_{i<j} g_{ij} -  10p.$$
Hence,  the Euler characteristic computation gives
$$ 1 - \beta_1(M^4) + \beta_2(M^4) -\beta_3(M^4) + 1 =   5 - \frac 1 3 \sum_{i<j<k} g_{ijk} + \frac 1 3 p .$$
Now, if $M^4$ is assumed to be simply-connected, $ 6 + 3 \beta_2(M^4) = 15 + p - \sum_{i<j<k} g_{ijk}$  follows;  since $g_{ijk} \ge 1$ trivially holds, we have $ 3 \beta_2(M^4) \le p-1.$

So,
\begin{equation}\label{eq:gem_complexity-beta2}
k(M^4) \ge 3 \beta_2 (M^4)
\end{equation}
may be stated, for each simply-connected PL 4-manifold $M^4$.

Now,  the classical theorems of Freedman and Donaldson~(\cite{Freedman-Quinn}) about the TOP classification of simply-connected closed 4-manifolds, together with more recent results by Furuta~(\cite{Furuta}),  ensure that intersection forms of type
$$ \pm 2n E_8 \oplus s \begin{pmatrix} 0 & 1 \\ 1 & 0 \end{pmatrix}$$
do represent PL $4$-manifolds  only if  $s > 2n $; hence, only PL 4-manifolds with  $\beta_2 (M^4)\ge 22$
occur in this case.
The thesis directly follows from the fact that $k(M^4) \le 65$ implies  $\beta_2 (M^4) \le 21$; so, only intersection forms of the two simplest types are allowed:
$$ r[1] \oplus r^{\prime}[-1]  \ \ \ \ \ {\text or} \ \ \ \ \ s \begin{pmatrix} 0 & 1 \\ 1 & 0 \end{pmatrix}$$
where $r+r^\prime = \beta_2(M^4)$ or $s = \frac 1 2\beta_2(M^4).$
\hskip 260pt\qed

As a consequence of the above result, we can already deduce the complete TOP classification of all PL 4-manifolds represented in our cystallization catalogues, i.e. up to gem-complexity 9.

We subdivide the results into two different statements, with respect to orientability assumptions;
in fact, in the non-orientable case, the PL classification of all PL 4-manifolds up to gem-complexity 9 follows, too (Proposition~\ref{Thm:classif_TOP_nonorientabile}), while in the orientable case the PL classification is possible only up to gem-complexity 5 (Proposition~\ref{Thm:classif_TOP_orientabile}):

\begin{proposition}  \label{Thm:classif_TOP_orientabile}
Let $M^4$ be an orientable PL 4-manifold. Then:
\ \par
\begin{itemize}
\item{} $k(M^4)=0\ \Longleftrightarrow\ M^4$ is PL-homeomorphic to $\mathbb S^4;$
\item{} $k(M^4)=3\ \Longleftrightarrow\ M^4$ is PL-homeomorphic to $\mathbb{CP}^{2};$
\item{} $k(M^4)=4\ \Longleftrightarrow\ M^4$ is PL-homeomorphic to $\mathbb{S}^{1} \times \mathbb{S}^{3};$
\item{} no orientable PL $4$-manifold $M^4$ exists  with $k(M^4) \in \{ 1,2,5 \};$
\item{} $k(\mathbb{CP}^{2}\#(\mathbb{S}^{1} \times \mathbb{S}^{3}))=7$, $k(\#_2(\mathbb{S}^{1} \times \mathbb{S}^{3}))=8$ and no other PL 4-manifold with handles $M^4$ exists with $k(M^4) \in \{6,7,8,9\}.$
\end{itemize}
Moreover, if $M^4$ is assumed to be handle-free, then:
\ \par
 \begin{itemize}
\item{} $k(M^4) \in \{6,7,8\} \ \Longrightarrow\ M^4$ is TOP-homeomorphic to either $\mathbb{S}^{2} \times \mathbb{S}^{2}$ or $\mathbb{CP}^{2} \# \mathbb{CP}^{2}$ or $\mathbb{CP}^{2} \# (- \mathbb{CP}^{2})$ or  $\mathbb{CP}^{2};$
\item{} $k(M^4)=9\ \Longrightarrow\ M^4$ is TOP-homeomorphic to either $\mathbb{S}^{2} \times \mathbb{S}^{2}$ or $(\#_r\mathbb {CP}^2)\#$ $(\#_{r^\prime}(-\mathbb {CP}^2))$, where $0 \leq r+r^\prime \leq 3.$
 \end{itemize}
 \end{proposition}

\dimo
In virtue of Proposition~\ref{Thm:rigid&dipole-free}, the gem-complexity of a closed connected orientable PL 4-manifold $M^4$ is congruent $mod \ 4$ to $p-1$, $2p$ being the order of a rigid dipole-free bipartite crystallization (i.e. an element of $\mathcal C^{(2p)}$). Moreover, by Proposition~\ref{Thm:calcolo_Eulero}, all elements of  $\bigcup_{1\leq p\leq 10}\mathcal C^{(2p)}$ represent simply-connected PL $4$-manifolds, which are obviously handle-free.

\par \noindent
Hence (as far as the PL category is concerned):
\begin{itemize}
\item[-] the statements for $k(M^4) \leq 3$ are direct consequences of the generation algorithm output, for $p \leq 4$ (see Proposition~\ref{Thm:calcolo_Eulero});
\item[-] $k(\mathbb{S}^{1} \times \mathbb{S}^{3})=4$ (resp. $k(\mathbb{CP}^{2}\#(\mathbb{S}^{1} \times \mathbb{S}^{3}))=7$) (resp. $k(\#_2(\mathbb{S}^{1} \times \mathbb{S}^{3}))=8$) follows from Proposition~\ref{Thm:rigid&dipole-free}(b), with $N^4= \mathbb S^4$, $h=1$ and $s=0$ (resp. $N^4= \mathbb{CP}^{2}$, $h=1$ and $s=0$) (resp.  $N^4= \mathbb S^4$, $h=2$ and $s=0$), since $k(\mathbb S^4)=0$, $k(\mathbb{CP}^{2})=3$ and no element of  $\bigcup_{1\leq p\leq 9}\mathcal C^{(2p)}$ represents a PL $4$-manifold with handles;
\item[-] for $k(M^4) \in \{4,5\}$, the statements follow from $k(\mathbb{S}^{1} \times \mathbb{S}^{3})=4$ and from the fact that $\bigcup_{5\leq p\leq 6}\mathcal C^{(2p)} = \emptyset$ (see Table~1);
\item[-]  the statement regarding the non-existence of other PL $4$-manifolds with handles (different from $\mathbb{CP}^{2}\#(\mathbb{S}^{1} \times \mathbb{S}^{3})$ and $\#_2(\mathbb{S}^{1} \times \mathbb{S}^{3})$) with  $k(M^4) = k \in \{6,7,8,9\}$ follows from Proposition~\ref{Thm:rigid&dipole-free}(b), too, together with the previous analysis concerning gem-complexity $k-4$ and $k-8$.
\end{itemize}

\smallskip

\par \noindent Finally, the statements involving TOP-homeomorphism follow immediately from Proposition~\ref{Thm:calcolo_Eulero} and Proposition~\ref{Thm:classif_TOP via gem-complexity}.
\hskip 320pt\qed

\begin{proposition}  \label{Thm:classif_TOP_nonorientabile}
Let $M^4$ be a non-orientable PL 4-manifold. Then:
\ \par
\begin{itemize}
\item{} $k(M^4)=4\ \Longleftrightarrow\ M^4$ is PL-homeomorphic to $\mathbb{S}^{1} \widetilde \times \mathbb{S}^{3};$
\item{} $k(M^4)=7\ \Longleftrightarrow\ M^4$ is PL-homeomorphic to either $\mathbb{RP}^{4}$ or $\mathbb{CP}^{2}\#(\mathbb{S}^{1} \widetilde \times \mathbb{S}^{3});$
\item{} $k(M^4)=8\ \Longleftrightarrow\ M^4$ is PL-homeomorphic to $\#_2(\mathbb{S}^{1} \widetilde \times \mathbb{S}^{3}).$
\end{itemize}
\noindent
Moreover, no non-orientable PL $4$-manifold $M^4$ exists with $k(M^4) \in \{0,1,2,3,5,6,9\}$.
\end{proposition}

\dimo
By Proposition~\ref{Thm:rigid&dipole-free}, the gem-complexity of a closed connected non-orientable PL 4-manifold $M^4$ is congruent $mod\ 4$ to $p-1$, $2p$ being the order of a rigid dipole-free crystallization (i.e. an element of $\mathcal C^{(2p)} \cup \mathcal{\tilde C}^{(2p)}$). Moreover, the generation algorithm output (Table~1) and Proposition~\ref{Thm:calcolo_Eulero} ensure that all elements of  $\bigcup_{1\leq p\leq 10} \left(\mathcal C^{(2p)} \cup \mathcal{\tilde C}^{(2p)}\right)$, except the standard order 16 crystallization of $\mathbb{RP}^{4}$, represent simply-connected PL $4$-manifolds; hence, no PL $4$-manifold with handles appears.

Now, the arguments are exactly the same as in the proof of Proposition~\ref{Thm:classif_TOP_orientabile}, when taking into account also the results of the generation algorithm for the non-bipartite case.
\hskip 40pt\qed

Another interesting consequence of the quoted results by Freedman, Donaldson and Furuta~(\cite{Freedman-Quinn}, \cite{Furuta}) is related to the topological classification of simply-connected PL 4-manifolds with respect to the invariant {\it regular genus}.

First, recall that the {\it genus} of a bipartite\footnote{An analogous definition exists in the non-bipartite case, too (see \cite{[G]});
for the purpose of the present work, however, the attention may be restricted to bipartite graphs.} $(n+1)$-coloured graph $(\Gamma, \gamma)$ with respect to a cyclic permutation $\e= (\e_o, \e_1, \dots, \e_{n-1}, \e_n=n)$ of $\Delta_n$ is the genus $\rho_\e (\Gamma)$ of the surface $F_\e$ into which $\Gamma$ regularly embeds (see \cite{[G]} for details); moreover,  $\rho_\e (\Gamma)$ may be directly computed by the following formula:
\begin{equation}\label{eq:genere}
\sum_{i\in \mathbb Z_{n+1}} g_{\e_{i}\e_{i+1}} + (1-n) \cdot p = 2 - 2 \rho_\e(\Gamma).
\end{equation}

Then, the {\it regular genus} of $\Gamma$  is defined as $\rho(\Gamma) = \min_\e \{\rho_\e(\Gamma)\}$, while the {\it regular genus} of an
orientable  PL $n$-manifold $M^n$ is defined as:
$$\mathcal G(M^n) = \min \{\rho(\Gamma) \ | \ (\Gamma, \gamma) \ \text{crystallization of} \  M^n\}.$$

\smallskip

With the above notations, the following statement holds:

\begin{proposition} \label{Thm:classif_TOP via regular-genus}
Any simply-connected PL 4-manifold $M^4$, with $\mathcal G(M^4)\leq 43$, is TOP-homeo\-morphic to
$$(\#_r\mathbb {CP}^2)\#(\#_{r^\prime}(-\mathbb {CP}^2))\quad or\quad\#_s(\mathbb S^2\times \mathbb S^2),$$
\noindent where $r+r^\prime = \beta_2(M^4),\ s = \frac 1
2\beta_2(M^4)$ and $\beta_2(M^4)$ is the second Betti number of $M^4$.
\end{proposition}

\dimo
Let $\Gamma$ be a crystallization of a PL 4-manifold $M^4$, and let $\rho_\e$ (resp. $ \rho_{\hat{\e_i}}$) be the regular genus of  $\Gamma$ (resp. $\Gamma_{\hat{\e_i}}$) with respect to a given cyclic permutation $\e$ of $\Delta_4$.
By applying formula (\ref{eq:genere}) both to $\Gamma$ and to each 4-residue of $\Gamma$, and by making use of formula (\ref{eq:Euler}), too, the following relations easily follow (see, for example, formulae (3)-(13) in the proofs of Lemma 1 and Lemma 2 in \cite{[C$_1992$]}, where they are given in the more general setting of crystallizations of bounded PL 4-manifolds):
\begin{itemize}
\item[(a)] $ g_{\e_{i-1},\e_i,\e_{i+1}}= 1 + \rho_\e -  \rho_{\hat{\e_{i+1}}} - \rho_{\hat{\e_{i+3}}}  \ \ \forall i\in \mathbb Z_4;$
\item[(b)] $\sum_{i\in \mathbb Z_4} g_{\e_{i-1},\e_i,\e_{i+1}} = 5 + 5 \rho_\e - 2 \sum_{i\in \mathbb Z_4}\rho_{\hat{\e_i}};$
\item[(c)] $  \chi (M^4) = 2 - 2 \rho_\e + \sum_{i\in \mathbb Z_4}\rho_{\hat{\e_i}}. $
\end{itemize}
Since $g_{\e_{i-1},\e_i,\e_{i+1}} \ge 1$ trivially holds,
inequality $2 \sum_{i\in \mathbb Z_4}\rho_{\hat{\e_i}} \le 5\rho_\e$ directly follows from (b).
By substituting it into (c), with the additional hypothesis $\pi_1(M^4)=0,$ we have $ 2 + \beta_2 (M^4) = \chi (M^4) \le 2 + [\frac {\rho_\e} 2] $
(where $[x]$ means the integer part of $x$), from which $ \beta_2 (M^4) \le [ \frac {\rho_\e} 2] $ follows.
As a consequence, relation
\begin{equation}\label{eq:gem_complexity-regular_genus}
\beta_2(M^4) \le \left[ \frac {\mathcal G(M^4)}2 \right]
\end{equation}
(already obtained in \cite[Proposition~2]{[Cav$_1992$]}, too) is proved to hold.

Now, the thesis directly follows from the fact that $\mathcal G(M^4)\leq 43$ implies  $\beta_2 (M^4) \le 21$, and from the already quoted well-known results about the classification of simply-connected closed 4-manifolds (exactly as in the proof of Proposition~\ref{Thm:classif_TOP via gem-complexity}).
\hskip 180pt\qed

As already pointed out, Proposition~\ref{Thm:classif_TOP via gem-complexity} and Proposition~\ref{Thm:classif_TOP via regular-genus}, together, prove Theorem~\ref{Thm:2} of Section~\ref{introduction}, related to the TOP classification of the PL 4-manifolds represented by crystallizations.

The (more interesting!) PL classification will be independently achieved in Section~\ref{sec:PL classification}, by making use of an implementation of the classifying algorithm described in Section~\ref{sec:paragrafo algoritmo classificazione}.

\section{Classification algorithm and possible applications}\label{sec:paragrafo algoritmo classificazione}

In order to complete the PL classification of the manifolds appearing in our catalogues of crystallizations, we exploit in the $n$-dimensional setting an idea already applied in dimension three.

\smallskip
First of all, let us call {\it admissible} a sequence of combinatorial moves which transforms a rigid dipole-free crystallization of a PL $n$-manifold $M$, into a rigid dipole-free crystallization of a PL $n$-manifold $M^\prime$ such that $M\cong_{PL}M^\prime\#_h \mathbb H$ ($h\ge 0$).

\medskip

Let now $X$ be a list of rigid dipole-free crystallizations; for any given set $\mathcal S$ of admissible sequences, it is possible to subdivide $X$ into equivalence classes with regard to $\mathcal S$.

More precisely, for each $\Gamma\in X$ and for each $\epsilon \in \mathcal S$, let $\theta_{\epsilon}(\Gamma)$ denote the (rigid dipole-free) crystallization obtained from $\Gamma$ by applying the admissible sequence $\epsilon$ of moves, and let us define the class of $\Gamma\in X$ with respect to $\mathcal S$ as:
$$cl_{\mathcal S}(\Gamma)=\{\Gamma^{\prime} \in X \ | \ \exists \epsilon, \epsilon^{\prime} \in \mathcal S, \ \theta_{\epsilon}(\Gamma) \ \text {and} \ \theta_{\epsilon^{\prime}}(\Gamma^{\prime}) \ \text {have the same code} \}$$

\smallskip

The following statement is a direct consequence of the above definition, together with Proposition~\ref{Thm:rho-pairs}:

\begin{proposition}
Given $\Gamma,\Gamma^{\prime}\in X$, if $cl_{\mathcal S}(\Gamma)=cl_{\mathcal S}(\Gamma^\prime)$, then there exist $h,k\in\mathbb N\cup\{0\}$ such that $|K(\Gamma)|\cong_{PL}M\#_h \mathbb H$ and $|K(\Gamma^\prime)|\cong_{PL}M\#_k \mathbb H$ (where the handles are orientable or not according to the bipartition of $\Gamma$ and $\Gamma^\prime$).
\end{proposition}

Note that no theoretical proof exists ensuring  that
$\ |K(\Gamma)|\cong_{PL}|K(\Gamma^\prime)|\ $ implies $cl_{\mathcal S}(\Gamma)=cl_{\mathcal S}(\Gamma^\prime)$
(as well as its generalization:
$|K(\Gamma)|\cong_{PL}|K(\Gamma^\prime)|\#_h \mathbb H \Rightarrow cl_{\mathcal S}(\Gamma)=cl_{\mathcal S}(\Gamma^\prime)$).

Nevertheless, in dimension three, existence has been proven of a set of admissible moves which are sufficient
to perform the topological (=PL) classification of all 3-manifolds admitting a coloured triangulation with at most $30$ tetrahedra
(\cite{[CC$_2$]}, \cite{[BCrG$_1$]}).

As we will see in the next section, the same turns out to be true,
in the 4-dimensional setting, for all elements of $\bigcup_{1\leq p\leq 9}\mathcal C^{(2p)}$, though with respect to a different set of moves.

In fact, the 3-dimensional classification algorithm employs dipole moves and $\rho$-pairs switchings (which are available in any dimension), together with {\it generalized dipole moves} (\cite{[FG]}), which are defined only for $n=3$.

Instead, in the $n$-dimensional setting we make use of a further set of moves, introduced by Lins and Mulazzani in \cite{[LM]}.

\begin{definition} Let $\Gamma$ be a gem of a PL $n$-manifold $M^n.$
Then:
\begin{itemize}
\item{} A {\it blob} is the insertion or cancellation of an $n$-dipole.
\item{} A {\it t-flip} is the switching of a pair $(e,f)$ of equally coloured edges which are both incident to an $h$-dipole $\Xi$ ($1\le h\le n-1$). An {\it s-flip} is the inverse move, i.e. the switching of a pair $(e,f)$ of equally coloured edges where either $e$ or $f$ belong to an $h$-dipole, which becomes an $(h-1)$-dipole after the transformation.  A {\it flip}
    is either an s- or a t-flip.
\end{itemize}
\end{definition}

\begin{figure}[th]
\begin{center}\scalebox{0.45}{\includegraphics{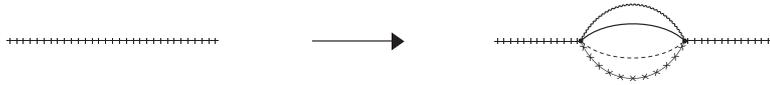}}\end{center}
\caption{\label{fig:blob_move} blob move}
\end{figure}

\bigskip

\begin{figure}[th]
\begin{center}\scalebox{0.45}{\includegraphics{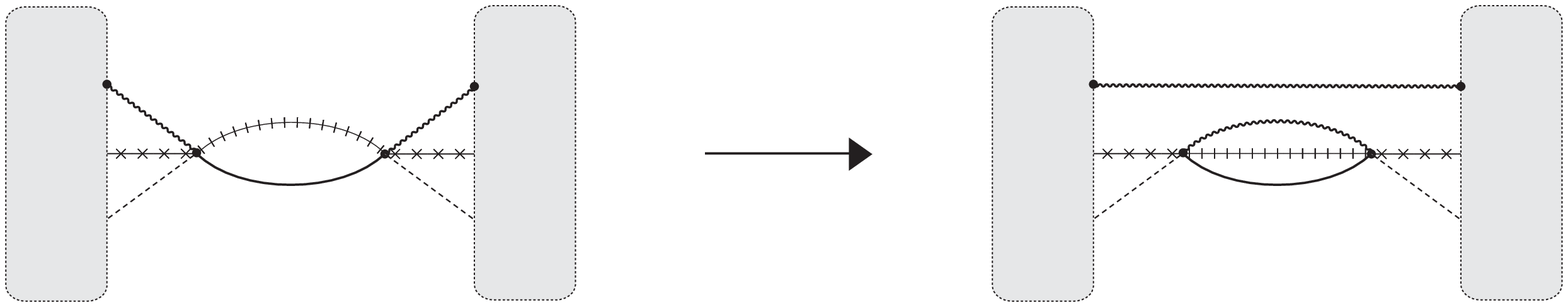}}\end{center}
\caption{\label{fig:flip_move} flip move}
\end{figure}

\bigskip

Flips and blobs on a gem do not change the represented manifold: \cite[Proposition~3]{[LM]}.
Actually, even if two crystallizations are known to represent the same manifold, there is no algorithmic procedure to determine a sequence of blobs and flips connecting them, nor an upper bound to the number of moves to be performed.
However, in the next section we will show that such an algorithm exists for the crystallizations appearing in our catalogues and that only one blob is sufficient.

\bigskip

In order to define the set of admissible moves $\bar{\mathcal S}$ which have been chosen to work on the catalogue $\bigcup_{1\leq p\leq 9}\mathcal C^{(2p)}$, let us introduce some definitions and notations.

Given an order $2p$ $(n+1)-$coloured graph $\Gamma$ there is a natural ordering of
its vertices induced by the \textit{rooted numbering algorithm}
generating its code (see \cite{[CG$_2$]}); so we can write
$V(\Gamma)=\{v_1,\ldots,v_{2p}\}$.

If $\Gamma$ is a rigid dipole-free crystallization of a PL $n$-manifold, given
$i\in\mathbb N_{2p}=\{1,\ldots,2p\}$,\ \ $c\in\Delta_n$, an $n$-tuple ${\bf x} = (x_1,\ldots,x_n)$ with $x_i\in\mathbb N_{2p}$ and a permutation $\tau$ of
$\hat c= \Delta_n - \{c\}$,
we denote by $\theta_{i,c,{\bf x},\tau}(\Gamma)$ the rigid dipole-free crystallization obtained from $\Gamma$ in the following way:
\begin{itemize}
\item [-] insert a blob over the $c$-coloured edge incident with $v_i$;

\item [-] for each $k\in\hat c$, consider, if exists, the s-flip on the pair of $\tau(k)$-coloured edges $(e,f)$, where $e$ belongs to the blob and $f$ is incident to $v_{x_k}$; then
perform the sequence of all possible s-flips of this type for increasing values of the parameter $k$;

\item [-] cancel dipoles and switch $\rho$-pairs in the resulting graph.
\end{itemize}

\medskip

$\theta_{i,c,{\bf x},\tau}$ obviously defines an admissible sequence.

\medskip

We denote by $\bar{\mathcal S}$ the set of all sequences $\theta_{i,c,{\bf x},\tau}$, where
$i\in\mathbb N_{2p},\ c\in\Delta_n,\ {\bf x}$ is an $n$-tuple of elements of $\mathbb N_{2p}$ and $\tau$ is a permutation of $\hat c$.

\begin{remark} Note that the above definition of $\bar{\mathcal S}$, as well as the classification algorithm itself, are independent from dimension. As a consequence, the partition into equivalence classes with respect to $\bar{\mathcal S}$ can be performed on any list of crystallizations of $n$-manifolds, in order to prove their PL-equivalence.
Furthermore, note that the PL manifold represented by an equivalence class is completely identified once at least one of the crystallizations of the class is ``known"; hence the algorithm can be also effective for the PL recognition of the manifolds involved in the list.
\end{remark}

\begin{remark} It is not difficult to prove that $\rho_n$-pairs cannot appear when the classification algorithm with respect to $\bar{\mathcal S}$ is applied to a set of bipartite crystallizations representing simply-connected $n$-manifolds (see Proposition~\ref{Thm:rho-pairs}(b)). Hence, in this case, $cl_{\mathcal S}(\Gamma)=cl_{\mathcal S}(\Gamma^\prime)$ surely implies $|K(\Gamma)|\cong_{PL}|K(\Gamma^\prime)|$.
 \end{remark}

In order to obtain PL classification results, the classification algorithm, with respect to $\bar{\mathcal S}$ and for $n=4$, has been implemented in a C++ program - called ``{\it $\Gamma 4$-class}".

As already pointed out in Section~\ref{introduction}, the program {\it $\Gamma 4$-class} can be applied to attempt to prove PL-equivalence between different (pseudo-)triangulations of the same topological 4-manifold; in fact, it is very easy to produce automatically a rigid dipole-free crystallization starting from any (pseudo-)triangulation (see Remark~\ref{rem:baricentrica} and Proposition~\ref{Thm:rigid&dipole-free}).

In particular, an application of {\it $\Gamma 4$-class} to the case of the $16$-vertices (resp. $17$-vertices) triangulation of the $K3$-surface (obtained in \cite{[CK]}  and \cite{[SK]} respectively) is in progress. The idea is similar to the one described in \cite{[BaS]}, \cite{[BuS]} and \cite{[BL]}, but the elementary moves involved in the automatic procedures are different ({\it blob} and {\it flips}, together with {\it dipole eliminations} and {\it $\rho$-pair switching}, instead of {\it edge-contraction} and {\it bistellar moves}). Hence, it is possible that one sequence succeeds when the others fail, or viceversa, with equal computational time employed.

\section{Classification results in PL=DIFF category}\label{sec:PL classification}

The application of the program {\it $\Gamma4$-class} to the catalogue $\bigcup_{1\leq p\leq 9}\mathcal C^{(2p)}$ yields the complete PL classification of the involved crystallizations (or, equivalently, of the dual contracted triangulations) as shown in the following proposition.\footnote{The proof of Proposition~\ref{Thm:biezione} shows that the PL classification performed through {\it $\Gamma4$-class} does not rely on the already stated results of Proposition~\ref{Thm:classif_TOP_orientabile}.}

\begin{proposition} \label{Thm:biezione} There is a bijective correspondence between the partition obtained by the program {\rm $\Gamma4$-class} and the set of PL 4-manifolds represented by $\bigcup_{1\leq p\leq 9} \left(\mathcal C^{(2p)} \cup \mathcal{\tilde C}^{(2p)}\right)$.
\noindent Moreover, the PL classification coincides with the TOP classification.
\end{proposition}

\dimo
Since there is only one non-bipartite crystallization in all our catalogues, the statement is trivial for $\bigcup_{1\leq p\leq 9}\mathcal{\tilde C}^{(2p)}.$

Program {\it $\Gamma4$-class} applied to the set $\bigcup_{1\leq p\leq 9}\mathcal C^{(2p)}$ produces a partition into five classes, which coincides with the partition induced by the second Betti number.
More precisely, all crystallizations with $\beta_2=0$ (resp. $\beta_2=1$) belong to the same class as the standard crystallization of $\mathbb S^4$ (resp. $\mathbb {CP}^2$), while the crystallizations with $\beta_2=2$ are subdivided into three classes, containing the standard crystallization of $\mathbb{CP}^{2} \# \mathbb{CP}^{2}$, $\mathbb{CP}^{2} \# (- \mathbb{CP}^{2})$ and $\mathbb{S}^{2} \times \mathbb{S}^{2}$ respectively.
\qed

The following proposition - whose statement already appeared in a partial and preliminary version  in \cite{oberwolfach} and in \cite{[CC_EN]} -  summarizes the complete PL classification of orientable (resp. non-orientable) PL 4-manifolds having gem-complexity up to 8 (resp. up to 9).

\begin{proposition} \label{Thm:PL-classificazione}
Let $M^4$ be a PL 4-manifold. Then:
\ \par
\begin{itemize}
\item{} $k(M^4)=0\ \Longleftrightarrow\ M^4$ is PL-homeomorphic to $\mathbb S^4;$
\item{} $k(M^4)=3\ \Longleftrightarrow\ M^4$ is PL-homeomorphic to $\mathbb{CP}^{2};$
\item{} $k(M^4)=4\ \Longleftrightarrow\ M^4$ is PL-homeomorphic to either $\mathbb{S}^{1} \times \mathbb{S}^{3}$ or $\mathbb{S}^{1} \widetilde \times \mathbb{S}^{3};$
\item{} $k(M^4)=6\ \Longleftrightarrow\ M^4$ is PL-homeomorphic to either $\mathbb{S}^{2} \times \mathbb{S}^{2}$ or $\mathbb{CP}^{2} \# \mathbb{CP}^{2}$ or $\mathbb{CP}^{2} \# (- \mathbb{CP}^{2});$
\item{} $k(M^4)=7\ \Longleftrightarrow\ M^4$ is PL-homeomorphic to either $\mathbb{RP}^{4}$  or $\mathbb{CP}^{2}\#(\mathbb{S}^{1} \times \mathbb{S}^{3})$ or $\mathbb{CP}^{2}\#(\mathbb{S}^{1} \widetilde \times \mathbb{S}^{3});$
\item{} $k(M^4)=8\ \Longleftrightarrow\ M^4$ is PL-homeomorphic to either $\#_2(\mathbb{S}^{1} \times \mathbb{S}^{3})$ or $\#_2(\mathbb{S}^{1} \widetilde \times \mathbb{S}^{3}).$
        \end{itemize}
Moreover:
\ \par
\begin{itemize}
\item{} no PL $4$-manifold $M^4$ exists with $k(M^4) \in \{1,2,5\};$
\item{} no exotic PL $4$-manifold exists, with $k(M^4) \le 8;$
\item{} any PL $4$-manifold $M^4$ with $k(M^4)=9$ is simply-connected (with second Betti number $\beta_2 \le 3$).
    \end{itemize}
\end{proposition}

\dimo
The statements concerning PL $4$-manifolds up to gem-complexity 8 are consequences of the previous Proposition~\ref{Thm:biezione}, together with Proposition~\ref{Thm:rigid&dipole-free}. The last statement, concerning gem-complexity 9, directly follows from Proposition~\ref{Thm:calcolo_Eulero}.
\hskip 180pt\qed

Note that the above Proposition~\ref{Thm:PL-classificazione} implies Theorem~\ref{Thm:1} (stated in Section~\ref{introduction}), when the attention is restricted to the handle-free PL $4$-manifolds.

\medskip

\begin{figure}[th]
\begin{center}\scalebox{0.7}{\includegraphics{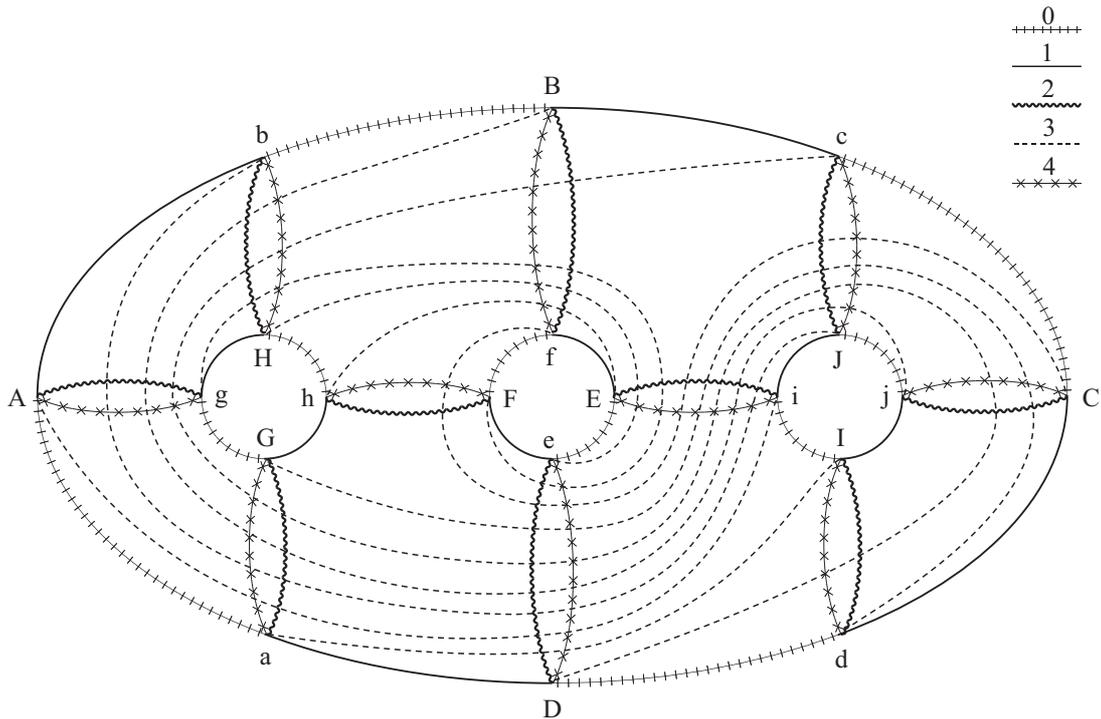}}\end{center}
\caption{\label{fig:S^4_20vertici} The order $20$ rigid dipole-free crystallization of $\mathbb S^4$}
\end{figure}

As a consequence of the (partial) analysis of the 4-dimensional crystallization catalogue $\bigcup_{1\leq p\leq 10}\mathcal C^{(2p)}$, together with a suitable application of the classification program {\it $\Gamma4$-class}, the following result may also be stated:

\begin{proposition}
A rigid crystallization of $\mathbb S^4$ exists, with $20$ vertices (see Figure~\ref{fig:S^4_20vertici}).
Apart from the standard order-two crystallization, it is the only rigid dipole-free crystallization of $\mathbb S^4$ up to $20$ vertices.
\end{proposition}

\dimo
As already stated in Proposition~\ref{Thm:calcolo_Eulero}, there exists exactly one crystallization $(\bar \Gamma, \bar \gamma) \in \bigcup_{7\leq p\leq 10}\mathcal C^{(2p)}$, with $\pi_1(|K(\bar \Gamma)|)=0$ and $\beta_2(|K(\bar \Gamma)|)=0$. Moreover, $\bar \Gamma$ has order $20$, while  no element of $ \bigcup_{2\leq p\leq 6}\mathcal C^{(2p)}$ represents a PL $4$-manifold with $\beta_2=0.$  Finally, by applying the program {\it $\Gamma4$-class} to the crystallization $\bar \Gamma$, the standard order two crystallization of $\mathbb S^4$ is obtained; hence, $|K(\bar \Gamma)| \cong_{PL} \mathbb S^4$ follows.
\hskip 350pt\qed

In \cite{[BaS]} the notion of {\it simple crystallization} of a (simply-connected) PL $4$-manifold is introduced, which is equivalent to the assumption $g_{rst}=1$ for any distinct $ r,s,t \in \Delta_4;$ moreover, the existence of simple crystallizations for any ``standard" simply-connected PL $4$-manifold (i.e.  $\mathbb S^4$, $\mathbb{CP}^{2}$, $\mathbb{S}^{2} \times \mathbb{S}^{2}$ and the $K3$-surface, together with the connected sums of them and/or their copies with opposite orientation) is proved.

As a consequence of our catalogues, we can state:

\begin{proposition}
\ \par
\begin{itemize}
\item{}  $\mathbb S^4$ and $\mathbb{CP}^{2}$ admit a unique simple crystallization;
\item{}  $\mathbb{S}^{2} \times \mathbb{S}^{2}$ admits exactly $267$ simple crystallizations;
\item{}  $\mathbb{CP}^{2} \# \mathbb{CP}^{2}$ admits exactly $583$ simple crystallizations;
\item{}  $\mathbb{CP}^{2} \# (- \mathbb{CP}^{2})$ admits exactly $258$ simple crystallizations.
    \end{itemize}
\rightline{$\Box$ \ \ }
\end{proposition}

Moreover, known facts and open problems about exotic structures on ``standard" simply-connected PL $4$-manifolds have the following implications on the existence of simple crystallizations:

\begin{proposition} \label{Thm:simple}
\ \par
\begin{itemize}
\item[(a)] Let $M^4$ be  $\mathbb S^4$ or $\mathbb{CP}^{2}$ or $\mathbb{S}^{2} \times \mathbb{S}^{2}$ or $\mathbb{CP}^{2} \# \mathbb{CP}^{2}$ or $\mathbb{CP}^{2} \# (- \mathbb{CP}^{2})$; if an exotic PL-structure on $M^4$ exists, then the corresponding PL-manifold does not admit a simple crystallization.
\item[(b)] Let $\ \bar M\ $ be a PL\ \ $4$-manifold\ \ TOP-homeomorphic\ \ but not\ \ PL-homeomorphic to\ \ $\mathbb{CP}^{2}\#_2(-\mathbb{CP}^{2})$; then, either $\bar M$ does not admit a simple crystallization, or $\bar M$ admits an order $20$ simple crystallization (i.e.: $k(\bar M)= 9 = k (\mathbb{CP}^{2} \#_2 (-\mathbb{CP}^{2}))$).
\item[(c)]  Let $r \in \{3,5,7,9,11,13 \} \cup \{r =4n-1 \ | \ n \ge 4\} \cup \{r =4n-2 \ | \ n \ge 23\}$; then, infinitely many simply-connected PL $4$-manifolds with $\beta_2=r$ do not admit a simple crystallization.
\end{itemize}
\end{proposition}

\dimo
Statements (a) and (b) directly follow from the fact - proved in \cite{[CCG]} - that simply-connected PL $4$-manifolds admitting a simple crystallization are characterized by $ k(M^4)= 3 \beta_2 (M^4).$

In order to prove statement (c), it is sufficient to recall the existence of infinitely many exotic structures TOP-homeomorphic to $\mathbb{CP}^{2} \#_4 (-\mathbb{CP}^{2})$, $\#_3\mathbb{CP}^2 \#_k (-\mathbb{CP}^2)$ for $k \in \{6, 8, 10\}$ and to $\#_{2n - 1}\mathbb{CP}^2 \#_{2n} (-\mathbb{CP}^2)$ for any integer $n \ge 1$ (see  \cite{[A-DP$_2010$]}), as well as to $\#_{2n-1} \mathbb{CP}^2 \#_{2n-1} (-\mathbb{CP}^2)$ for $n \ge 23$ (see \cite{[A-I-DP$_2013$]}).  Then, the statement follows from the above characterization of simply-connected PL $4$-manifolds admitting simple crystallizations (due to  \cite{[CCG]}), together with the (obvious) finiteness of PL $4$-manifolds having a fixed gem-complexity.
\hskip 110pt\qed

Proposition~\ref{Thm:simple}(b) suggests possible interesting consequences of the work in progress on the catalogue $\mathcal C^{(20)}$: in fact, the PL-characterization of all PL $4$-manifolds represented by order $20$ crystallizations could yield results about exotic structures on simply-connected 4-manifolds with $\beta_2 \leq 3$ and/or about the existence of simple crystallizations, in case $\beta_2 =3$.

More generally, we hope that further developments in the generation and classification of 4-dimensional crystallization catalogues, for increasing gem-complexity, could be useful to face open problems concerning different PL-structures on the same TOP $4$-manifold.

\subsection*{Acknowledgements} We thank the referee for his/her helpful comments.
We acknowledge the CINECA award under the ISCRA initiative, for the availability of high performance computing resources and support.

\end{document}